\documentclass[11pt]{article}
\usepackage{a4}
\usepackage[fleqn]{amsmath}
\usepackage{amsthm}
\usepackage{amssymb}

\newtheorem{proposition}{Proposition}
\newtheorem{corollary}[proposition]{Corollary}
\newtheorem{theorem}[proposition]{Theorem}
\newtheorem{lemma}[proposition]{Lemma}

\theoremstyle{definition}
\newtheorem{definition}{Definition}
\newtheorem{example}{Example}
\newtheorem{remark}{Remark}

\newcommand{\goth}[1]{\mathfrak{#1}}
\newcommand{\B}[1]{\mathbf{#1}}

\newcommand{\Dynr}{\ensuremath{\mathbf{Dynr}}}
\newcommand{\Homsp}{\ensuremath{\mathbf{Homsp}}}
\newcommand{\Map}{\ensuremath{\mathbf{Map}}}

\newcommand{\leftinvariant}[1]{\ensuremath{\overrightarrow{#1}}}
\newcommand{\rightinvariant}[1]{\ensuremath{\overleftarrow{#1}}}

\newcommand{\ovr}{\overrightarrow}

\newcommand{\g}{\mathfrak g}
\newcommand{\h}{\mathfrak h}
\newcommand{\n}{\mathfrak{n}}
\newcommand{\bb}{\mathfrak{b}}
\renewcommand{\u}{\mathfrak u}
\renewcommand{\v}{\mathfrak v}
\renewcommand{\sl}{\mathfrak{sl}}

 \DeclareMathOperator{\Ad}{Ad}
\DeclareMathOperator{\id}{id} \DeclareMathOperator{\Alt}{Alt}
\DeclareMathOperator{\Lie}{Lie} 
 \DeclareMathOperator{\CYB}{CYB}
\DeclareMathOperator{\ad}{ad} \DeclareMathOperator{\Tr}{Tr}
\DeclareMathOperator{\Ker}{Ker} \DeclareMathOperator{\gr}{gr}
\DeclareMathOperator{\End}{End} \DeclareMathOperator{\spanv}{span}
\DeclareMathOperator{\arccoth}{arccoth}
\DeclareMathOperator{\Dom}{Dom} \DeclareMathOperator{\Fun}{Fun}
\DeclareMathOperator{\Hom}{Hom} \DeclareMathOperator{\fin}{fin}
\DeclareMathOperator{\rrr}{r} \DeclareMathOperator{\diag}{diag}
\DeclareMathOperator{\Aut}{Aut}

\def\ftext#1{{\let\thefootnote\relax
\footnotetext{\vskip-.9\baselineskip\noindent #1}}}

\begin{document}
\title{Dynamical Yang-Baxter equations, quasi-Poisson homogeneous spaces,
and quantization}
\author{E.~Karolinsky$^{\,\diamond}$, K.~Muzykin, A.~Stolin,
and V.~Tarasov$^{\,*}$}
\date{}

\maketitle

\ftext{$^{\diamond}$Supported in part by the Royal Swedish Academy of Sciences
\\
$^*$Supported in part by RFFI grant 02--01--00085a and
CRDF grant RM1--2334--MO--02}

\begin{sloppy}

\vskip-.6\baselineskip
\hrule height0pt
\thispagestyle{empty}

\section{Introduction}

This paper is a continuation of \cite{EKASrmatr}. Let us recall
the main result of \cite{EKASrmatr}. Let $G$ be a Lie group,
$\g=\Lie G$, $U\subset G$ a connected closed Lie subgroup such
that the corresponding subalgebra $\u\subset\g$ is reductive in
$\g$ (i.e., there exists an $\u$-invariant subspace
$\goth{m}\subset\g$ such that $\g=\u\oplus\goth{m}$), and
$\Omega\in(\goth{u}\otimes\goth{u})\oplus(\goth{m}\otimes\goth{m})$
a symmetric tensor. Take a solution $\rho\in\g\otimes\g$ of the
classical Yang-Baxter equation such that $\rho+\rho^{21}=\Omega$
and consider the corresponding Poisson Lie group structure
$\pi_\rho$ on $G$. Assuming additionally that
\begin{equation}\label{old_cond}
\rho+s\in\frac{\Omega}{2}+\left(\textstyle{\bigwedge}^2\goth{m}\right)^\u
\end{equation}
for some element $s\in\bigwedge^2\g$ that satisfies a certain
``twist'' equation, we establish a 1-1 correspondence between the
moduli space of classical dynamical r-matrices for the pair
$(\g,\u)$ with the symmetric part $\frac{\Omega}{2}$ and the set
of all structures of Poisson homogeneous $(G,\pi_\rho)$-space on
$G/U$. We emphasize that the first example of such a
correspondence was found by Lu in \cite{LuCDrmatr}.

We develop the results of \cite{EKASrmatr} principally in two
directions. First, we generalize the main result of
\cite{EKASrmatr}.
We replace Poisson Lie groups (resp.\ Poisson homogeneous spaces)
by quasi-Poisson Lie groups (resp.\ quasi-Poisson homogeneous
spaces), but even in the Poisson case our result (see Theorem
\ref{rmatr_bijection_theorem}) is stronger than in
\cite{EKASrmatr}: condition \eqref{old_cond} is relaxed now. We
hope that now we present this result in its natural generality.

Secondly, we propose a partial quantization of the results of
\cite{EKASrmatr}. We explain how, starting from dynamical twist
for a pair $(U\g,\h)$ (where $\g$ is a Lie algebra, $\h$ is its
abelian subalgebra, and $U\g$ is the universal enveloping algebra
of $\g$), one can get a $G$-equivariant star-product on $G/H$
(where $H\subset G$ are connected Lie groups corresponding to
$\h\subset\g$). In the case $\g$ is complex simple and $\h$ is its
Cartan subalgebra we give a representation-theoretic explanation
of our results in terms of Verma modules. We also provide an
analogue of these results for quantum universal enveloping
algebras. Notice that results in this direction were obtained in
the recent papers \cite{DM} and \cite{AL}. However, our approach
is completely elementary (cf.\ \cite{DM}), and we emphasize the
connection between star-products and dynamical twists very
explicitly (cf.\ \cite{AL}). We also propose a method that allows
one (under certain conditions) to obtain non-dynamical twists from
dynamical ones. This is a quantization of the ``classical'' result
obtained in \cite[Appendix B]{EKASrmatr}.

Let us explain the structure of this paper in more details.
Section \ref{classical} is devoted to the quasi-classical picture.
In Subsection \ref{classical_main} we remind the definitions of
classical dynamical r-matrices, quasi-Poisson Lie groups and their
quasi-Poisson homogeneous spaces, and then formulate and prove the
main result of this section, Theorem
\ref{rmatr_bijection_theorem}. In Subsection
\ref{semisimple_subsection} we consider an example: the case of
quasi-triangular (in the strict sense) classical dynamical
r-matrices for the pair $(\g,\u)$, where $\g$ is a complex
semisimple Lie algebra, and $\u$ is its regular reductive
subalgebra. Section \ref{star} contains the construction of
star-products from dynamical twists. As an example, we write down
an explicit formula for an equivariant star-product on a regular
semisimple coadjoint orbit of $SL(2)$ (note that this formula was
also obtained as an example in \cite{AL}; similar formulas
appeared earlier in physical papers) and observe its relation to
certain Verma modules. In Section \ref{Verma} we give a more
conceptual explanation of the connection between equivariant
quantization and Verma modules; our approach differs from that of
\cite{AL}. In Section \ref{q-Verma} we obtain an analogue of the
results of Section \ref{Verma} for the case of quantum universal
enveloping algebras, therefore providing some examples of
``quantum homogeneous spaces'' related to dynamical twists for
QUE-algebras.
Finally, in Section \ref{ndtwists} we present a way
from dynamical twists to ``usual'' (non-dynamical) ones; we apply
this construction to an explicit calculation of the universal
twist for the universal enveloping algebra of two-dimensional
nonabelian Lie algebra.

All Lie algebras in this paper assumed to be finite dimensional,
and the ground field (except of Section \ref{q-Verma}) is $\mathbb
C$.

When the paper was finished we got sad news about unexpected
passing away of Joseph Donin. We dedicate this paper to his
memory.

\subsection*{Acknowledgements} The authors are grateful to Maria
Gorelik, Aleksandr Il'inskii, and Leonid Vaksman for useful
discussions on the topic of the paper. Part of this research was
done during E.K.'s visit to University of G\"oteborg, the
University of Hong Kong, Peking University, and E.K.'s and K.M.'s
visits to Free University of Berlin; we thank our colleagues there
for their hospitality.

\section{Quasi-Poisson homogeneous spaces and classical dynamical
r-matrices}\label{classical}

\subsection{General results}\label{classical_main}

In this section we describe a connection between quasi-Poisson
homogeneous spaces and classical dynamical r-matrices (see Theorem
\ref{rmatr_bijection_theorem}).

First we recall some definitions. Suppose $G$ is a Lie group,
$U\subset G$ its connected Lie subgroup. Let $\goth{g}$ and
$\goth{u}$ be the corresponding Lie algebras. Choose a basis $x_1,
\ldots, x_r$ in $\goth{u}$. Denote by $D$ the formal neighborhood
of zero in $\goth{u}^*$. By functions from $D$ to a vector space
$V$ we mean elements of the space $V[[x_1, \ldots, x_r]]$, where
$x_i$ are regarded as coordinates on $D$. Further, if $\omega \in
\Omega^k(D, V)$ is a $k$-form on $D$ with values in vector space
$V$, then by $\overline{\omega}: D \rightarrow \bigwedge^k\goth{u}
\otimes V$ we denote the corresponding function.

\begin{definition}[see \cite{ES_moduli_spaces}]
\emph{Classical dynamical r-matrix for the pair $(\goth{g},
\goth{u})$} is an $\goth{u}$-equivariant function $r: D
\rightarrow \goth{g}\otimes\goth{g}$ that satisfies the
\emph{classical dynamical Yang-Baxter equation (CDYBE)}:
\[
\Alt(\overline{dr})+\CYB(r)=0,
\]
where $\CYB(r)=[r^{12},r^{13}]+[r^{12},r^{23}]+[r^{13},r^{23}]$,
and for $x \in \g^{\otimes 3}$ we set $\Alt(x) = x^{123} + x^{231}
+ x^{312}$.
\end{definition}

We will also require the \emph{quasi-unitarity property}:
\[
r + r^{21} = \Omega \in (S^2\goth{g})^\goth{g}.
\]
It is easy to see that if $r$ satisfies the CDYBE and the
quasi-unitarity condition, then $\Omega$ is constant.

We denote the set of all classical dynamical r-matrices for the
pair $(\goth{g}, \goth{u})$ such that $r + r^{21} = \Omega$ by
$\Dynr(\goth{g}, \goth{u}, \Omega)$.

Denote by $\Map(D,G)^\goth{u}$ the set of all
$\goth{u}$-equivariant maps from $D$ to $G$. Suppose that $r: D
\rightarrow \goth{g\otimes\goth{g}}$ is an $\goth{u}$-equivariant
function. Then for any $g\in\Map(D,G)^\goth{u}$ define a function
$r^g: D \rightarrow \goth{g}\otimes\goth{g}$ by
\[
r^g= (\Ad_g \otimes \Ad_g)(r-
\overline{\eta_g}+\overline{\eta_g}^{21}+\tau_g),
\]
where $\eta_g = g^{-1}dg$, and $\tau_g(\lambda) = (\lambda \otimes
1 \otimes 1) ([ \overline{\eta_g}^{12} ,
\overline{\eta_g}^{13}](\lambda))$. Then $r^g$ is a classical
dynamical r-matrix if and only if $r$ is. The transformation $r
\mapsto r^g$ is called a {\it gauge transformation}. In fact, it
is an action of the group $\Map(D,G)^\goth{u}$ on $\Dynr(\goth{g},
\goth{u}, \Omega)$.

Following \cite{ES_moduli_spaces}, we denote the moduli space
$\Map_0(D,G)^\goth{u} \backslash \Dynr(\goth{g}, \goth{u},
\Omega)$ by $\mathcal{M}(\goth{g}, \goth{u}, \Omega)$ (here
$\Map_0(D,G)^\goth{u} = \left\{ g\in\Map(D,G)^\goth{u}: g(0)=e
\right\}$).

Now we recall the definition of quasi-Poisson Lie groups and their
quasi-Poisson homogeneous spaces (for details see
\cite{YKS_qbialg_qpoiss, AAYKS_mpairs_mmaps, AAYKSEM}).

\begin{definition}
Let $G$ be a Lie group, $\goth{g}$ its Lie algebra, $\pi_G$ a
bivector field on $G$, and $\varphi\in\bigwedge^3\goth{g}$. A
triple $(G,\pi_G,\varphi)$ is called a \emph{quasi-Poisson Lie
group} if
\begin{gather*}
\pi_G(gg')=(l_g)_*\pi_G(g')+(r_{g'})_*\pi_G(g),\\
\frac{1}{2}[\pi_G,\pi_G]=\rightinvariant{\varphi}-\leftinvariant{\varphi},\\
[\pi_G,\rightinvariant{\varphi}]=0,
\end{gather*}
where $l_g$ (resp. $r_g$) is left (resp. right) multiplication by
$g$, $\leftinvariant{a}$ (resp. $\rightinvariant{a}$) is the left
(resp. right) invariant tensor field on $G$ corresponding to $a$
and $[.,.]$ is the Schouten bracket of multivector fields.
\end{definition}

\begin{definition}
Suppose that $(G,\pi_G,\varphi)$ is a quasi-Poisson group, $X$ is
a homogeneous $G$-space equipped with a bivector field $\pi_X$.
Then $(X,\pi_X)$ is called a \emph{quasi-Poisson homogeneous
$(G,\pi_G,\varphi)$-space} if
\begin{gather*}
\pi_X(gx)=(l_g)_*\pi_X(x)+(\rho_x)_*\pi_G(g),\\
\frac{1}{2}[\pi_X,\pi_X]=\varphi_X
\end{gather*}
(here $l_g$ denotes the mapping $x\mapsto g\cdot x$, $\rho_x$ is
the mapping $g\mapsto g\cdot x$, and $\varphi_X$ is the trivector
field on $X$ induced by $\varphi$).
\end{definition}

Now take $\rho \in \goth{g}\otimes\goth{g}$ such that
$\rho+\rho^{21}=\Omega\in(S^2\g)^\g$. Let $\Lambda = \rho-
\frac{\Omega}{2} \in \bigwedge^2\goth{g}$. Define a bivector field
on $G$ by $\pi_\rho = \leftinvariant{\rho} - \rightinvariant{\rho}
= \leftinvariant{\Lambda} - \rightinvariant{\Lambda}$. Set
$\varphi = \varphi_\rho = - \CYB(\rho)$. Then
$(G,\pi_\rho,\varphi)$ is a quasi-Poisson Lie group (such
quasi-Poisson Lie groups are called \emph{quasi-triangular}).
Denote by $\Homsp(G,\pi_\rho,\varphi,U)$ the set of all
$(G,\pi_\rho,\varphi)$-homogeneous quasi-Poisson structures on
$G/U$. We will see that, under certain conditions, there is a
bijection between $\mathcal{M}(\goth{g}, \goth{u}, \Omega)$ and
$\Homsp(G,\pi_\rho,\varphi,U)$.

Assume that $b \in (\goth{g}\otimes\goth{g})^\goth{u}$ is such
that $b + b^{21} = \Omega$. Let $B = b - \frac{\Omega}{2}$. Define
a bivector field on $G$ by $\tilde{\pi}_b^\rho =
\leftinvariant{b}-\rightinvariant{\rho} =
\leftinvariant{B}-\rightinvariant{\Lambda}$. Then there is a
bivector field on $G/U$ defined by $\pi_b^\rho(\underline{g}) =
p_*(\tilde{\pi}_b^\rho(g))$ (here $p:G\to G/U$ is the canonical
projection, and $\underline{g}=p(g)$). It is well defined, since
$b$ is $\goth{u}$-invariant.

\begin{proposition}\label{CYB0_lemma} 
In this setting $(G/U, \pi_b^\rho)$ is a $(G, \pi_\rho,
\varphi)$-quasi-Poisson homogeneous space iff $\CYB(b) = 0$ in
$\bigwedge^3(\goth{g}/\goth{u})$.
\end{proposition}
\begin{proof}
First we check the ``multiplicativity'' of $\pi_b^\rho$. For all
$g \in G, u\in U$ we have
\begin{equation*}
g\cdot\tilde{\pi}_b^\rho(u) + \pi_\rho(g)\cdot u = gu\cdot b -
\rho\cdot gu = \tilde{\pi}_b^\rho(gu).
\end{equation*}
Using
$p_*$, we get the required equality $\pi_b^\rho(\underline{g}) =
g\cdot \pi_b^\rho(\underline{e}) + p_*\pi_\rho(g)$.

Now we need to prove that $\frac{1}{2}[\pi_b^\rho,\pi_b^\rho] =
\varphi_{G/U}$ iff $\CYB(b) = 0$ in
$\bigwedge^3(\goth{g}/\goth{u})$. We check this directly:
\begin{gather*}
\frac{1}{2}[\tilde{\pi}_b^\rho,\tilde{\pi}_b^\rho] =
\frac{1}{2}\left( [\leftinvariant{B},\leftinvariant{B}] +
[\rightinvariant{\Lambda},\rightinvariant{\Lambda}] \right) =
-\leftinvariant{\CYB(B)}+\rightinvariant{\CYB(\Lambda)} =\\
-\leftinvariant{\CYB(b)}+\rightinvariant{\varphi}.
\end{gather*}
Consequently, $\frac{1}{2}[\pi_b^\rho,\pi_b^\rho] =
p_*(-\leftinvariant{\CYB(b)}+\rightinvariant{\varphi}) =
-p_*(\leftinvariant{\CYB(b)}) + \varphi_{G/U}$. So we see that
$\frac{1}{2}[\pi_b^\rho,\pi_b^\rho] = \varphi_{G/U}$ iff $\CYB(b)
= 0$ in $\bigwedge^3(\goth{g}/\goth{u})$.
\end{proof}

Suppose $r \in \Dynr(\goth{g},\goth{u},\Omega)$.

\begin{proposition}[see \cite{LuCDrmatr}]
$\CYB(r(0)) = 0$ in $\bigwedge^3(\goth{g}/\goth{u})$.\qed
\end{proposition}

\begin{corollary}
The correspondence $r\mapsto \pi_{r(0)}^\rho$ is a map from
$\Dynr(\goth{g},\goth{u},\Omega)$ to
$\Homsp(G,\pi_\rho,\varphi,U)$.\qed
\end{corollary}

\begin{proposition}[see \cite{EKASrmatr}]
If $g\in \Map_0(D,G)^\goth{u}$, then $\pi_{r(0)}^\rho =
\pi_{r^g(0)}^\rho$.\qed
\end{proposition}

\begin{corollary}
The correspondence $r\mapsto \pi_{r(0)}^\rho$ defines a map from
$\mathcal{M}(\goth{g},\goth{u},\Omega)$ to
$\Homsp(G,\pi_\rho,\varphi,U)$.\qed
\end{corollary}

Now consider the following conditions:
\begin{subequations}
\begin{gather}\label{complement_a}
\text{$\goth{u}$ has an $\goth{u}$-invariant complement $\goth{m}$
in $\goth{g}$};\\
\label{complement_b} \Omega \in
(\goth{u}\otimes\goth{u})\oplus(\goth{m}\otimes\goth{m}).
\end{gather}
\end{subequations}

Assume that \eqref{complement_a} holds and consider the algebraic
variety
\begin{equation*}
\mathcal{M}_\Omega = \left\{ x\in \frac{\Omega}{2} +
\left(\textstyle{\bigwedge}^2\goth{m}\right)^\goth{u}\,\bigg|\,\CYB(x)
= 0\text{ in }\textstyle{\bigwedge}^3(\goth{g}/\goth{u})\right\}.
\end{equation*}

\begin{theorem}[Etingof, Schiffman; see \cite{ES_moduli_spaces}]\label{Etingof_Schiffman_rmatr_theorem}
$(1)$ Any class $\mathcal{C} \in
\mathcal{M}(\goth{g},\goth{u},\Omega)$ has a representative $r \in
\mathcal{C}$ such that $r(0) \in \mathcal{M}_\Omega$. Moreover,
this defines an embedding $\mathcal{M}(\goth{g},\goth{u},\Omega)
\rightarrow \mathcal{M}_\Omega$.

$(2)$ Assume that \eqref{complement_b} holds. Then the map
$\mathcal{M}(\goth{g},\goth{u},\Omega) \rightarrow
\mathcal{M}_\Omega$ defined above is a bijection.\qed
\end{theorem}

\begin{proposition} \label{MOmega_Homsp_bijection_proposition}
Assume that \eqref{complement_a} holds. Then the mapping $b\mapsto
\pi_b^\rho$ from $\mathcal{M}_\Omega$ to
$\Homsp(G,\pi_\rho,\varphi,U)$ is a bijection.
\end{proposition}
\begin{proof}
Let's build the inverse mapping. Assume that $\pi$ is a bivector
field on $G/U$ defining a structure of a
$(G,\pi_\rho,\varphi)$-quasi-Poisson homogeneous space. Then
$\pi(\underline{e}) \in \bigwedge^2(\goth{g}/\goth{u}) =
\bigwedge^2\goth{m}$. Consider
$b=\frac{\Omega}{2}+\pi(\underline{e})+p_*(\Lambda)$. We will
prove that $b\in\mathcal{M}_\Omega$ and the mapping $\pi\mapsto b$
is inverse to the mapping $b\mapsto\pi_b^\rho$.

First we prove that $b\in
(\bigwedge^2\goth{m})^\goth{u}+\frac{\Omega}{2}$. For all $u\in U$
we have $\pi(\underline{e}) + p_*(\Lambda) =
\pi(u\cdot\underline{e}) + p_*(\Lambda\cdot u) =
u\cdot\pi(\underline{e}) + p_*(\pi_\rho(u)) + p_*(\Lambda\cdot u)
= u\cdot\pi(\underline{e}) +
p_*(u\cdot\rho-u\cdot\frac{\Omega}{2}) = u\cdot(\pi(\underline{e})
+ p_*(\Lambda))$. This means that $\pi(\underline{e}) +
p_*(\Lambda) \in (\bigwedge^2\goth{m})^\goth{u}$.

Now we prove that $\pi = \pi_b^\rho$. By definition,
$\pi_b^\rho(\underline{g}) = p_*( g\cdot\pi(\underline{e}) +
g\cdot p_*\Lambda - \Lambda\cdot g) =
\pi(\underline{g})+p_*(g\cdot p_*\Lambda-\Lambda\cdot
g-g\cdot\Lambda+\Lambda\cdot g) = \pi(\underline{g})$. So
$\pi_b^\rho$ defines a structure of
$(G,\pi_\rho,\varphi)$-quasi-Poisson homogeneous space. By
Proposition \ref{CYB0_lemma}, this means that
$b\in\mathcal{M}_\Omega$.
\end{proof}

\begin{theorem}\label{rmatr_bijection_theorem}
Suppose \eqref{complement_a} and \eqref{complement_b} are
satisfied. Then the map $r\mapsto \pi_{r(0)}^\rho$ from
$\mathcal{M}(\goth{g},\goth{u},\Omega)$ to
$\Homsp(G,\pi_\rho,\varphi,U)$ is a bijection.
\end{theorem}
\begin{proof}
This theorem follows from Theorem
\ref{Etingof_Schiffman_rmatr_theorem} and Proposition
\ref{MOmega_Homsp_bijection_proposition}.
\end{proof}

\begin{remark}
If $\varphi = -\CYB(\rho) = 0$, then $(G, \pi_\rho)$ is a Poisson
Lie group. In this case we get a bijection between
$\mathcal{M}(\goth{g},\goth{u},\Omega)$ and the set of all Poisson
$(G,\pi_\rho)$-homogeneous structures on $G/U$.
\end{remark}

\begin{remark}
Assume that only \eqref{complement_a} holds. Clearly, in this case
the map $r\mapsto \pi_{r(0)}^\rho$ defines an embedding
$\mathcal{M}(\goth{g},\goth{u},\Omega)\hookrightarrow\Homsp(G,\pi_\rho,\varphi,U)$.
\end{remark}

\begin{remark}
If \eqref{complement_a} fails, then the space
$\mathcal{M}(\goth{g},\goth{u},\Omega)$ may be
infinite-dimensional (see \cite{Xu_triang_rmatr}), while
$\Homsp(G,\pi_\rho,\varphi,U)$ is always finite-dimensional.
\end{remark}

\subsection{Example: the semisimple case}\label{semisimple_subsection}
Assume that $\goth{g}$ is a semisimple Lie algebra. Choose a
Cartan subalgebra $\goth{h}\subset\goth{g}$ and denote by $\B R$
the corresponding root system. Suppose $\left<.,.\right>$ is a
non-degenerate symmetric invariant bilinear form on $\goth{g}$,
and $\Omega\in (S^2\goth{g})^\g$ is the corresponding tensor. We
will describe $\mathcal{M}_\Omega$ for a reductive Lie subalgebra
$\goth{u}\subset\goth{g}$ containing $\goth{h}$.

Precisely, consider a set $\B U\subset \B R$ such that $\goth{u} =
\goth{h} \oplus \sum_{\alpha \in \B U} \goth{g}_\alpha$ is a
reductive Lie subalgebra. In this case we will call $\B U$
\emph{reductive} (in other words, a set $\B U \subset\B R$ is
reductive iff $(\B U+\B U)\cap \B R\subset \B U$ and $-\B U=\B
U$). Note that in this situation condition \eqref{complement_a} is
satisfied, since $\goth{m} = \sum_{\alpha \in \B R \backslash \B
U}\goth{g}_\alpha$ is an $\goth{u}$-invariant complement to
$\goth{u}$ in $\goth{g}$.

Fix $E_\alpha \in \goth{g}_\alpha$ such that $\left< E_\alpha ,
E_{-\alpha} \right> = 1$ for all $\alpha \in \B R$. Then $\Omega =
\Omega_\goth{h} + \sum_{\alpha \in \B R} E_\alpha \otimes
E_{-\alpha}$, where $\Omega_\goth{h} \in S^2\goth{h}$. Notice that
\eqref{complement_b} is also satisfied.

\begin{proposition}
Suppose that $x=\sum_{\alpha\in \B R}x_\alpha E_\alpha \otimes
E_{-\alpha}$. Then $x + \frac{\Omega}{2} \in \mathcal{M}_\Omega$
iff
\begin{subequations}
\begin{gather}
x_\alpha=0 \text{ for }\alpha \in \B U;
\label{mOmega_a}\\
x_{-\alpha}=-x_\alpha \text{ for } \alpha \in \B R;
\label{mOmega_b}\\
\text{if }\alpha,\beta\in \B R\backslash \B U, \gamma\in \B U, \alpha+\beta+\gamma=0,
\text{then }x_\alpha+x_\beta=0;
\label{mOmega_c}\\
\text{if }\alpha,\beta,\gamma\in \B R\backslash \B U, \alpha+\beta+\gamma=0,
\text{then } x_\alpha x_\beta + x_\beta x_\gamma + x_\gamma x_\alpha = -1/4.
\label{mOmega_d}
\end{gather}
\end{subequations}
\end{proposition}

Note that \eqref{mOmega_c} is equivalent to the following
condition:
\begin{equation*}
\text{if }\alpha\in \B R\backslash \B U, \beta\in \B U,\text{then
}x_{\alpha+\beta}=x_\alpha.
\end{equation*}

\begin{proof}
It is easy to see that $x \in (\bigwedge^2\goth{m})^\goth{h}$ iff
\eqref{mOmega_a} and \eqref{mOmega_b} are satisfied.

Suppose that $c_{\alpha\beta}$ are defined by $[E_\alpha, E_\beta]
= c_{\alpha\beta} E_{\alpha + \beta}$.

For any $\gamma \in \B U$ we have
\begin{gather*}
\left[ E_\gamma , x \right]= \sum_{\alpha \in \B R\backslash \B U}
x_\alpha( [E_\gamma,E_\alpha]\otimes E_{-\alpha} +
E_\alpha\otimes[E_\gamma,E_{-\alpha}] ) =\\
\sum_{\alpha,\beta \in \B R\backslash \B U, \alpha+\beta+\gamma=0}
\left( x_\alpha c_{\gamma\alpha}E_{-\beta}\otimes E_{-\alpha} -
x_\alpha c_{\gamma\alpha} E_{-\alpha}\otimes E_{-\beta} \right) = \\
\sum_{\alpha,\beta \in \B R\backslash \B U,
\alpha+\beta+\gamma=0}\left( x_\beta c_{\gamma\alpha} - x_\alpha
c_{\gamma\beta} \right)E_{-\alpha}\otimes E_{-\beta} =\\
\sum_{\alpha,\beta \in \B R\backslash \B U,
\alpha+\beta+\gamma=0}( x_\alpha + x_\beta ) c_{\gamma\alpha}
E_{-\alpha}\otimes E_{-\beta}.
\end{gather*}
Thus $x$ is $\goth{u}$-invariant if and only if $x_\alpha +
x_\beta = 0$ for all $\alpha, \beta \in \B R\backslash \B U$ such
that $\alpha + \beta\in \B U$.

Finally, we calculate $\CYB\left(x + \frac{\Omega}{2}\right) =
\CYB(x) + \CYB\left(\frac{\Omega}{2}\right)$ (see
\cite{AAYKS_mpairs_mmaps}):

\begin{gather*}
\CYB(x) =\\
\sum_{\alpha, \beta \in \B R} x_\alpha x_\beta
\bigl([E_\alpha,E_\beta]\otimes E_{-\alpha} \otimes E_{-\beta} +
E_\alpha \otimes [E_{-\alpha},E_\beta] \otimes E_{-\beta} +\\
E_\alpha \otimes E_\beta \otimes [E_{-\alpha},E_{-\beta}]\bigr) =\\
\sum_{\alpha, \beta, \gamma \in \B R, \alpha + \beta + \gamma = 0}
\bigl( x_\alpha x_\beta c_{\alpha\beta} E_{-\gamma}\otimes
E_{-\alpha}\otimes E_{-\beta}-\\
x_\alpha x_\beta c_{\alpha\beta} E_{-\alpha}\otimes
E_{-\gamma}\otimes E_{-\beta}+x_\alpha x_\beta c_{\alpha\beta}
E_{-\alpha}\otimes E_{-\beta}\otimes E_{-\gamma}
\bigr)=\\
\sum_{\alpha, \beta, \gamma \in \B R, \alpha + \beta + \gamma = 0}
c_{\alpha\beta}(x_\alpha x_\beta + x_\alpha x_\gamma + x_\beta
x_\gamma) E_{-\alpha} \otimes E_{-\beta} \otimes E_{-\gamma},\\ \\
\CYB\left(\frac{\Omega}{2}\right) \equiv \frac{1}{4}\sum_{\alpha,
\beta, \gamma \in \B R\backslash \B U, \alpha + \beta + \gamma =
0} c_{\alpha\beta}E_{-\alpha} \otimes E_{-\beta} \otimes
E_{-\gamma}
\\ \text{(mod $\goth{u}\otimes\goth{g}\otimes\goth{g} +
\goth{g}\otimes\goth{u}\otimes\goth{g} +
\goth{g}\otimes\goth{g}\otimes\goth{u}$)}.
\end{gather*}

So the image of $\CYB\left(x + \frac{\Omega}{2}\right)$ in
$\bigwedge^3(\goth{g}/\goth{u})$ vanishes if and only if the
condition \eqref{mOmega_d} is satisfied.
\end{proof}

\begin{proposition}\label{x_alpha_proposition}
Suppose $\Pi \subset \B R$ is a set of simple roots, $\B R_+$ is
the corresponding set of positive roots. Choose a subset $\Delta
\subset \Pi$ such that $\B N=(\spanv\Delta)\cap \B R$ contains $\B
U$. Find $h \in \mathfrak{h}$ such that $\alpha(h)\notin \pi i
\mathbb{Z}$ for $\alpha\in \B N\backslash \B U$ and $\alpha(h)\in
\pi i \mathbb{Z}$ for $\alpha\in \B U$. Then $x_\alpha$ defined by
\[
x_\alpha=\begin{cases}0, & \alpha\in \B U\\
\frac{1}{2}\coth\alpha(h),&\alpha\in \B N\backslash \B U\\
\pm 1/2, & \alpha\in \pm \B R_+\backslash \B N
\end{cases}
\]
satisfies \eqref{mOmega_a}--\eqref{mOmega_d}. Moreover, any
function satisfying \eqref{mOmega_a}--\eqref{mOmega_d} is of this
form.
\end{proposition}

First we prove the second part of the proposition. Set
\begin{equation*}
\B P=\{ \alpha\,|\, x_\alpha \neq -1/2 \}.
\end{equation*}
It is obvious that $\B U \subset \B P$.
\begin{lemma}
$\B P$ is parabolic.
\end{lemma}
\begin{proof}
Obviously, $\B P \cup (-\B P) = \B R$.

We have to prove that if $\alpha, \beta \in \B P$ and
$\alpha+\beta \in \B R$, then $\alpha+\beta \in \B P$. We do it by
considering several cases. If $\alpha, \beta \in \B U$, then
$\alpha+\beta \in \B U \subset \B P$. If $\alpha \in \B P
\backslash \B U$ and $\beta \in \B U$, then $x_{\alpha+\beta} =
x_\alpha \neq -1/2$ by \eqref{mOmega_c} and $\alpha+\beta \in \B
P$. If $\alpha, \beta \in \B P \backslash \B U$, there are two
possibilities. If $\alpha+\beta \in \B U$, then there is nothing
to prove. If $\alpha+\beta \notin \B U$, then, by
\eqref{mOmega_d}, $x_\alpha x_\beta - x_{\alpha+\beta}( x_\alpha +
x_\beta ) = -1/4$. If $x_{\alpha+\beta} = -1/2$, then from this
equation it follows that $x_\alpha = -1/2$. Consequently,
$\alpha+\beta \in \B P$.
\end{proof}

Since $\B P$ is parabolic, there exists a set of positive roots
$\Pi \subset \B R$ and a subset $\Delta \subset \Pi$ such that $\B
P = \B R_+ \cup \B N$ (see \cite{B_groupes_algebres}, chapter VI,
$\S$ 1, proposition 20); here $\B R_+$ is the set of positive
roots corresponding to $\Pi$, and $\B N = (\spanv\Delta) \cap \B
R$ is the Levi subset corresponding to $\Delta$.

%

Let $\B N_+ = \B N\cap \B R_+$ be the set of positive roots in $\B
N$ corresponding to $\Delta$. For all $\alpha \in \Delta\backslash
\B U$ let $y_\alpha=\arccoth 2 x_\alpha$, for $\alpha \in
\Delta\cap \B U$ let $y_\alpha=0$. Find $h \in \mathfrak{h}$ such
that $y_\alpha=\alpha(h)$. Now we prove that $h$ satisfies
Proposition \ref{x_alpha_proposition}.
\begin{lemma}
$\alpha(h)\notin \pi i \mathbb{Z}$ and $x_\alpha =
\frac{1}{2}\coth\alpha(h)$ for all $\alpha \in \B N\backslash \B
U$; $\alpha(h)\in \pi i \mathbb{Z}$ for $\alpha \in \B U$.
\end{lemma}
\begin{proof}
It is enough to prove this for $\alpha$ positive, so that we can
use the induction on the length $l(\alpha)$. The case
$l(\alpha)=1$ is trivial. Suppose that $l(\alpha)=k$. Then we can
find $\alpha'\in \B N_+$ and $\alpha_k \in \Delta$ such that
$l(\alpha')=k-1$ and $\alpha=\alpha'+\alpha_k$. Consider two
cases.

First, suppose that $\alpha \in \B U$.

If $\alpha_k \in \B U$, then $\alpha'\in \B U$. By induction,
$\alpha(h)=\alpha'(h) \in \pi i \mathbb{Z}$.

If $\alpha_k \notin \B U$, then $\alpha' \notin \B U$. By
induction assumption, $x_{\alpha'}=\frac{1}{2}\coth\alpha'(h)$.
From \eqref{mOmega_c} it follows that
$0=x_{\alpha'}+x_{\alpha_k}=\frac{1}{2}(\coth\alpha'(h)+\coth\alpha_k(h))$
and, consequently, $\alpha(h)\in \pi i \mathbb{Z}$.

Now suppose that $\alpha \notin \B U$.

If $\alpha_k \in \B U$, then $\alpha'\notin \B U$. Since
$\alpha_k(h)=0$, by \eqref{mOmega_c} we have
$x_\alpha=x_{\alpha'+\alpha_k}=x_{\alpha'}=
\frac{1}{2}\coth\alpha'(h)=\frac{1}{2}\coth\alpha(h)$.

When $\alpha_k \notin \B U$, then there are two possibilities
again. If $\alpha'\in \B U$, then by induction $\alpha'(h) \in \pi
i \mathbb{Z}$. By \eqref{mOmega_c}, $0=x_\alpha+x_{-\alpha_k}$.
Consequently,
$x_\alpha=x_{\alpha_k}=\frac{1}{2}\coth\alpha_k(h)=\frac{1}{2}\coth\alpha(h)$.
If $\alpha'\notin \B U$, then, by \eqref{mOmega_d}, $x_\alpha
x_{-\alpha'}+x_{-\alpha'}x_{-\alpha_k}+x_{-\alpha_k}x_\alpha=-1/4$.
This equation can be rewritten as
\begin{equation*}
x_\alpha=
\frac{1/4+x_{\alpha'}x_{\alpha_k}}{x_{\alpha'}+x_{\alpha_k}}=
\frac{1}{2}\cdot\frac{1+\coth\alpha'(h)\coth\alpha_k(h)}{\coth\alpha'(h)+\coth\alpha_k(h)}=
\frac{1}{2}\coth\alpha(h),
\end{equation*}
and the lemma is proved.
\end{proof}

To prove the first part of the proposition we need the following
root theory lemma.

\begin{lemma}\label{Y_lemma}
Suppose $\B P\subset \B R$ is parabolic. Then $\B Y=\B R\backslash
\B P$ has the following properties:
\begin{subequations}
\begin{gather}
(-\B Y)\cap \B Y=\emptyset;
\label{Y_a}\\
(\B Y+\B Y)\cap \B R \subset \B Y;
\label{Y_b}\\
\text{if } \alpha\in \B Y, \beta \in \B R\backslash \B Y
\text{ and } \alpha-\beta\in \B R,\text{ then } \alpha-\beta \in \B Y.
\label{Y_c}
\end{gather}
\end{subequations}
\end{lemma}
\begin{proof}
Since \eqref{Y_a} is obvious and \eqref{Y_b} follows from
\eqref{Y_a} and \eqref{Y_c}, we prove only the last property: if
$\alpha \in \B Y$ and $\beta \in \B P$ are such that $\alpha-\beta
\in \B P$, then, since $\B P$ is parabolic, we would have $\alpha
= (\alpha-\beta) + \beta \in \B P$. So $\alpha-\beta \in \B Y$.
\end{proof}

Now we just check \eqref{mOmega_a}--\eqref{mOmega_d} directly.
Suppose that $\B N$ is defined as in the proposition. Let $\B Y=\B
R_+\backslash \B N$. Then $\B P=\B R\backslash \B Y = -\B R_+ \cup
\B N$ is a parabolic set, and $\B Y$ satisfies
\eqref{Y_a}--\eqref{Y_c}.

\begin{lemma}
Suppose $x_\alpha$ is as defined in Proposition
\ref{x_alpha_proposition}. Then $x_\alpha$ satisfies
\eqref{mOmega_a}--\eqref{mOmega_d}.
\end{lemma}
\begin{proof}
\eqref{mOmega_a} we have already, \eqref{mOmega_b} is trivial.

To prove \eqref{mOmega_c}, consider the following cases. First,
take $\alpha, \beta \in \B N\backslash \B U$, $\gamma \in \B U$,
$\alpha+\beta+\gamma=0$. Then $x_\alpha+x_\beta =
\frac{1}{2}(\coth\alpha(h)+\coth\beta(h)) = 0$ as $\alpha+\beta =
-\gamma \in \B U$. The case $\alpha \in \B N\backslash \B U, \beta
\in \B R\backslash \B N, \gamma \in \B U$ is impossible, because
then we would have $\beta = -\alpha-\gamma \in \B N$. The case
$\alpha, \beta \in \pm \B Y, \gamma \in \B U$ is also impossible,
because $-\gamma = \alpha+\beta \in \pm \B Y$. Finally, if
$\alpha\in \pm \B Y, \beta\in \mp \B Y, \gamma \in \B U$, then
$x_\alpha+x_\beta=\pm \frac{1}{2} \mp \frac{1}{2}=0$.

Condition \eqref{mOmega_d} can be proved in a similar way.
%
%
\end{proof}

Now to summarize:

\begin{theorem}\label{semisimple_quasipoissclassif}
Suppose $U \subset G$ is the connected Lie subgroup corresponding
to $\goth{u} \subset \goth{g}$. Take $\rho \in
\goth{g}\otimes\goth{g}$ such that $\rho+\rho^{21}=\Omega$ and set
$\varphi = - \CYB(\rho)$. Then any
$(G,\pi_\rho,\varphi)$-homogeneous quasi-Poisson space structure
on $G/U$ is exactly of the form $\pi = \pi_{x+\Omega/2}^\rho$ for
some $x=\sum_{\alpha \in R} x_\alpha E_\alpha \otimes
E_{-\alpha}$, where $x_\alpha$ is defined in Proposition
\ref{x_alpha_proposition}. \qed\end{theorem}

\begin{remark}
Let $\rho$ be any solution of the classical Yang-Baxter equation
such that $\rho+\rho^{21}=\Omega$ (see \cite{BD_FAn}). Then
$(G,\pi_\rho)$ is a Poisson Lie group and therefore Theorem
\ref{semisimple_quasipoissclassif} provides the list of all
$(G,\pi_\rho)$-homogeneous Poisson space structures on $G/U$.
\end{remark}

\begin{remark}
In \cite{D_poiss_hom_spaces}, Drinfeld assigned to each point of
any Poisson homogeneous space a Lagrangian (i.e., maximal
isotropic) subalgebra in the corresponding double Lie algebra.
Roughly speaking, this construction gives a one-to-one
correspondence between Poisson homogeneous spaces up to
isomorphism and Lagrangian subalgebras up to conjugation. In fact,
literally the same is true in the quasi-Poisson case (see
\cite{EKKMquasi}). In our situation for any
$\rho\in\goth{g}\otimes\goth{g}$ such that $\rho+\rho^{21}=\Omega$
the Manin pair that corresponds to the (quasi-)Poisson Lie group
$(G,\pi_\rho,\varphi_\rho)$ is the same and equals
$(\g\times\g,\g_{\diag})$; here $\g\times\g$ is equipped with the
invariant scalar product
\[
Q((x,y),(x',y'))=\langle x,x'\rangle-\langle y,y'\rangle,
\]
and $\g_{\diag}$ is the image of the diagonal embedding
$\g\to\g\times\g$. By \cite{D_poiss_hom_spaces, EKKMquasi} we have
a bijection between $(G,\pi_\rho,\varphi_\rho)$-homogeneous
(quasi-)Poisson space structures on $G/U$ and
Lagrangian subalgebras $\goth{l}\subset\g\times\g$ such that
$\goth{l}\cap\g_{\diag}=\u_{\diag}$ (this $\goth{l}$ corresponds
by Drinfeld to the base point $\underline{e}\in G/U$). Since all
(quasi-)Poisson Lie groups $(G,\pi_\rho,\varphi_\rho)$ are related
by twisting, we conclude that the
Lagrangian subalgebra $\goth{l}$ corresponding to
$\pi_{x+\Omega/2}^\rho$ defined in Theorem
\ref{semisimple_quasipoissclassif} is independent of $\rho$ (see
\cite{EKKMquasi, EKASrmatr}). Using Drinfeld's definition, it is
easy to compute $\goth{l}$ by taking, for example,
$\rho=\frac{\Omega}{2}$:

\begin{proposition}
Under the notation of Proposition \ref{x_alpha_proposition} and
Theorem \ref{semisimple_quasipoissclassif} let
\[
\goth{b}_\pm=\h\oplus\sum_{\alpha\in\B R_+}\g_{\pm\alpha}, \ \
\n=\h\oplus\sum_{\alpha\in\B N}\g_{\alpha},
\]
$\goth{p}_\pm=\n+\goth{b}_\pm$, and
$\theta=\exp(2\ad_h)\in\Aut\n$. Then
\[
\goth{l}=\{(x,y)\in\goth{p}_-\times\goth{p}_+\,|\,\theta(p_-(x))=p_+(y)\},
\]
where $p_\pm:\goth{p}_\pm\to\n$ are the canonical projections.\qed
\end{proposition}

\end{remark}

\section{Dynamical twists and equivariant
quan\-ti\-za\-ti\-on}\label{star}

Let $\g$ be a Lie algebra, $\h$ its abelian subalgebra. We will
use the standard Hopf algebra structure on the universal
enveloping algebra $U\g$ (and denote by $\Delta$ the
comultiplication, by $\varepsilon$ the counit, etc.). Suppose a
meromorphic function
\begin{equation*}
J:\h^*\to (U\g\otimes U\g)^\h[[\hbar]]
\end{equation*}
is a {\it quantum dynamical twist}, i.e.,
\begin{equation}\label{qdtwist}
J(\lambda)^{12,3}J(\lambda-\hbar
h^{(3)})^{12}=J(\lambda)^{1,23}J(\lambda)^{23}
\end{equation}
and
\begin{equation}\label{qdtwist2}
(\varepsilon\otimes\id)(J(\lambda))=(\id\otimes\varepsilon)(J(\lambda))=1.
\end{equation}
Here
\begin{gather*}
J(\lambda-\hbar
h^{(3)})^{12}=\\
J(\lambda)\otimes1-\hbar\sum_i\frac{\partial
J}{\partial\lambda_i}(\lambda)\otimes
h_i+\frac{\hbar^2}{2}\sum_{i,j}\frac{\partial^2J}{\partial\lambda_i\partial\lambda_j}(\lambda)\otimes
h_ih_j-\ldots ,
\end{gather*}
$h_i$ form a basis in $\h$, and $\lambda_i$ are the corresponding
coordinates on $\h^*$. We also use the standard notation
$A^{12,3}=(\Delta\otimes\id)(A)$, $A^{23}=1\otimes A$, etc., for
any $A\in(U\g\otimes U\g)[[\hbar]]$.

Let $G$ be a connected Lie group that corresponds to $\g$. Assume
that there exists the closed connected subgroup $H\subset G$
corresponding to $\h$. Identify $C^\infty(G/H)$ with right
$H$-invariant smooth functions on $G$.

Fix any $\lambda\in\Dom J$ and for any $f_1,f_2\in C^\infty(G/H)$
define $f_1\star_\lambda
f_2=\overrightarrow{J(\lambda)}(f_1,f_2):=(m\circ\overrightarrow{J(\lambda)})(f_1\otimes
f_2)$, where $\overrightarrow{J(\lambda)}$ is the left-invariant
differential operator corresponding to $J(\lambda)\in(U\g\otimes
U\g)^\h[[\hbar]]$, and $m$ is the usual multiplication in
$C^\infty(G/H)$ (extended naturally on $C^\infty(G/H)[[\hbar]]$).
Since $J(\lambda)$ is $\h$-invariant, we have $f_1\star_\lambda
f_2\in C^\infty(G/H)[[\hbar]]$. As usual, we extend
$\star_\lambda$ on $C^\infty(G/H)[[\hbar]]$.

\begin{theorem}
The correspondence $(f_1,f_2)\mapsto f_1\star_\lambda f_2$ is a
$G$-equivariant star-product on $G/H$.
\end{theorem}

\begin{proof}
We have
\begin{gather*}
(f_1\star_\lambda f_2)\star_\lambda
f_3=\left(m\circ(m\otimes\id)\circ\overrightarrow{J(\lambda)^{12,3}J(\lambda)^{12}}\right)(f_1\otimes
f_2\otimes f_3),\\
f_1\star_\lambda(f_2\star_\lambda
f_3)=\left(m\circ(m\otimes\id)\circ\overrightarrow{J(\lambda)^{1,23}J(\lambda)^{23}}\right)(f_1\otimes
f_2\otimes f_3).
\end{gather*}
Since
\begin{equation*}
J(\lambda-\hbar h^{(3)})^{12}-J(\lambda)^{12}\in(U\g\otimes
U\g\otimes U\g\cdot\h)[[\hbar]],
\end{equation*}
we see that
\begin{equation*}
J(\lambda)^{12,3}J(\lambda)^{12}\equiv
J(\lambda)^{1,23}J(\lambda)^{23}\mod(U\g\otimes U\g\otimes
U\g\cdot\h)[[\hbar]],
\end{equation*}
and
$\overrightarrow{J(\lambda)^{12,3}J(\lambda)^{12}}=\overrightarrow{J(\lambda)^{1,23}J(\lambda)^{23}}$.
This proves the associativity of $\star_\lambda$. The conditions
$f\star_\lambda1=1\star_\lambda f=f$ follows from
(\ref{qdtwist2}). The $G$-equivariance of $\star_\lambda$ is
obvious.
\end{proof}

\begin{example}\label{sl2ex}
Suppose $\g=\sl(2)$, $G=SL(2)$. Let $x,y,h$ be the standard basis
in $\g$, and $\h=\mathbb Ch$.

Consider the ABRR quantum dynamical twist $J$ for $(\g,\h)$, i.e.,
\begin{equation}
J(\lambda)=1+\sum_{n\geq1}J_n(\lambda),
\end{equation}
where
\begin{gather*}
J_n(\lambda)=\\
\frac{(-1)^n}{n!}\hbar^ny^n\otimes(\lambda+\hbar(n+1-h))^{-1}
\ldots(\lambda+\hbar(2n-h))^{-1}x^n=\\
\frac{(-1)^n}{n!}\hbar^ny^n\otimes x^n(\lambda-\hbar
h)^{-1}(\lambda-\hbar(h+1))^{-1}\ldots(\lambda-\hbar (h+n-1))^{-1}
\end{gather*}
(see \cite{ABRR, ES_DYBE}). Here we identify $\lambda\in\h^*$ with
$\lambda(h)\in\mathbb C$. Notice that $J(\lambda)$ is defined for
$\lambda\neq0$.

Clearly, for any $f_1,f_2\in C^\infty(G/H)$ we have
\begin{equation}\label{starlambda}
f_1\star_\lambda
f_2=f_1f_2+\sum_{n\geq1}\overrightarrow{J_{\lambda,\hbar}^{(n)}}(f_1,f_2),
\end{equation}
where
\begin{equation}\label{starlambdacomponents}
J_{\lambda,\hbar}^{(n)}=\frac{(-1)^n}{n!}
\frac{\hbar^n}{\lambda(\lambda-\hbar)\ldots(\lambda-(n-1)\hbar)}y^n\otimes
x^n.
\end{equation}

Equip $\g$ with the invariant scalar product defined by $\langle
a,b\rangle=\Tr ab$. Let us identify $\g$ with $\g^*$ and $\h$ with
$\h^*$ via $\langle\cdot,\cdot\rangle$. One can also consider
$G/H$ as a (co)adjoint $G$-orbit $O_\lambda$ of
$\lambda\in\h^*\subset\g^*$ (or $\frac{\lambda}{2}h\in\g$). Denote
by $f_a$ the restriction onto $O_\lambda$ of the linear function
on $\g^*$ generated by $a\in\g$ (i.e., in terms of $G/H$ we have
$f_a(g)=\frac{\lambda}{2}\langle gHg^{-1},a\rangle$). It is clear
that $\overrightarrow{x^n}f_a=\overrightarrow{y^n}f_a=0$ for all
$a\in\g$ and $n\geq2$. Therefore
\begin{equation*}
f_a\star_\lambda f_b-f_b\star_\lambda f_a=\hbar
\overrightarrow{u_\lambda}(f_a,f_b)=\hbar f_{[a,b]},
\end{equation*}
where
\begin{equation*}
u_\lambda=\frac{1}{\lambda}(x\otimes y-y\otimes x).
\end{equation*}
In other words, the quasiclassical limit of $\star_\lambda$ is
exactly $O_\lambda$ equipped with the Kirillov-Kostant-Souriau
bracket.

Let us now restrict ourselves to regular functions on $O_\lambda$.
Note that for any two such functions $f_1,f_2$ the series
(\ref{starlambda}) has only finitely many non-vanishing summands.
This allows us to fix the ``deformation parameter'' in
(\ref{starlambda}) and (\ref{starlambdacomponents}) (i.e., set
formally $\hbar=1$). Of course, this makes sense if
$\lambda\not\in\mathbb Z_+$. Denote by $A_\lambda$ the obtained
algebra, i.e., the algebra of regular functions on $O_\lambda$
with the multiplication $\star_\lambda$ (and $\hbar=1$).

It is not hard to check directly that for any $a,b\in\g$ we have
\begin{equation}\label{fafb}
f_a\star_\lambda
f_b=\left(1-\frac{1}{\lambda}\right)f_af_b+\frac{1}{2}f_{[a,b]}
+\frac{\lambda}{2}\langle a,b\rangle.
\end{equation}
Iterating (\ref{fafb}), we see that $\star_\lambda$ is compatible
with the standard (i.e., by polynomial degree) filtration on
$A_\lambda$. Since
\begin{equation*} f_a\star_\lambda
f_b-f_b\star_\lambda f_a=f_{[a,b]},
\end{equation*}
we get an algebra homomorphism $F: U\g\to A_\lambda$ defined by
$a\mapsto f_a$ for all $a\in\g$. Obviously, $F$ is filtered with
respect to the standard filtrations on $U\g$ and $A_\lambda$.

Consider the Casimir element $c=xy+yx+\frac{1}{2}h^2\in U\g$. Let
us calculate $F(c)$. We have
\begin{gather*}
f_x\star_\lambda f_y+f_y\star_\lambda
f_x+\frac{1}{2}f_h\star_\lambda f_h=\\
\left(1-\frac{1}{\lambda}\right)\left(2f_xf_y+\frac{1}{2}f_h^2\right)+\frac{3\lambda}{2}=
\left(1-\frac{1}{\lambda}\right)\frac{\lambda^2}{2}+\frac{3\lambda}{2}=\frac{\lambda(\lambda+2)}{2}.
\end{gather*}
Therefore $c-\frac{\lambda(\lambda+2)}{2}\in\Ker F$, and $F$
induces the (filtered) homomorphism
$\widetilde{F}:U\g\Big/\!\left(c-\frac{\lambda(\lambda+2)}{2}\right)\to
A_\lambda$.

Let us now pass to the corresponding gradings. Since
$\gr\widetilde{F}$ is obviously surjective on each graded
component, and the dimensions of the corresponding graded
components of
$U\g\Big/\!\left(c-\frac{\lambda(\lambda+2)}{2}\right)$ and
$A_\lambda$ are the same, we see that $\gr\widetilde{F}$ is an
isomorphism. Thus $\widetilde{F}$ is also an isomorphism, i.e.,
\begin{equation*}
A_\lambda\simeq
U\g\bigg/\!\left(c-\frac{\lambda(\lambda+2)}{2}\right).
\end{equation*}

Notice that for $\lambda\not\in\mathbb Z_+$ the ideal
$\left(c-\frac{\lambda(\lambda+2)}{2}\right)$ is exactly the
kernel of the natural homomorphism $U\g\to\End M(\lambda)$, where
$M(\lambda)$ is the Verma module with highest weight $\lambda$.
Therefore we get an embedding $A_\lambda\hookrightarrow\End
M(\lambda)$.
\end{example}

\section{Verma modules and equivariant quantization}\label{Verma}

In this section we give an explanation of the appearance of Verma
modules in Example \ref{sl2ex}.

\subsection{General construction}

Let $F=\mathbb C[G]$ be the algebra of polynomial functions on a
simple complex Lie group $G$, which is an algebra generated by
matrix elements of finite dimensional representations. Set
$\g=\Lie G$. We equip $F$ by a structure of $U\g$-module algebra
via $(a,f)\mapsto\ovr{a}f$.



Let $H$ be a Cartan subgroup of $G$, and $\h=\Lie H$. We also
define $\Fun(G/H)=F[0]=\{f\in L\,|\, \ovr{h}f=0 \mbox{ for any }h
\in \h\}$.


Let $M$ be a $\g$-module. On $\mbox{Hom}_\g(M, M\otimes F)$ we
introduce a natural structure of algebra in the following way. Let
$\varphi, \psi \in \mbox{Hom}_\g(M, M\otimes F)$. Set
\begin{equation}\label{basic_star_prod}
\varphi \ast \psi =(\id\otimes m)\circ (\varphi \otimes \id)\circ
\psi,
\end{equation}
where $m$ is the multiplication in $F$. We notice that $\varphi
\ast \psi \in \mbox{Hom}_\g(M, M \otimes F)$. It is not difficult
to see that the multiplication $\ast$ is associative and
$e(m)=m\otimes 1$ is the identity element.

\begin{remark}
{\rm Any linear map $f:M \to M \otimes F$ can be considered as a
function $\widetilde{f}:G \to \End M$ in the following way:
$\widetilde{f}(g)(v)=f(v)(g)$. We notice that
$\widetilde{\varphi\ast
\psi}=\widetilde{\varphi}\cdot\widetilde{\psi}$. Elements of
$\mbox{Hom}_\g(M, M\otimes F)$ can be distinguished by the
following lemma.}
\end{remark}

\begin{lemma}\label{diffeqn}
Let $\varphi \in \Hom(M, M \otimes F)$. Then $\varphi \in
\Hom_\g(M, M \otimes F)$ if and only if the corresponding function
$\widetilde{\varphi}$ satisfies the first order differential
equation $\ovr{a}\widetilde{\varphi}(g)=[\widetilde{\varphi}(g),
a_M]$.\qed
\end{lemma}

\begin{corollary} There exists an embedding $\Hom_\g(M, M\otimes F)
\hookrightarrow \End M$ given by the formula $\varphi \mapsto
\widetilde{\varphi}(e)$.
\end{corollary}

\begin{proof} The fact that this map is a homomorphism follows from
the remark above. Injectivity follows from Lemma \ref{diffeqn}.
\end{proof}

\medskip

Let us identify $\Hom_\g(M, M\otimes F)$ with its image in $\End
M$ when appropriate.

\begin{proposition}
There exists a homomorphism $U\g \rightarrow \Hom_\g(M, M\otimes
F)$ given by the formula $x \mapsto \Ad_{g^{-1}}(x)_M$.
\end{proposition}

\begin{proof} We have to verify that $\Ad_{g^{-1}}(x)$ satisfies the
differential equation from Lemma \ref{diffeqn}, which is
straightforward.
\end{proof}

\begin{remark}
{\rm The composition
\[
U\g \rightarrow \Hom_\g(M, M\otimes F) \rightarrow \End M
\]
is the standard homomorphism $U\g\rightarrow \End M$, $x\mapsto
x_M$.}
\end{remark}

\subsection{Verma modules}

Now fix a triangular decomposition $\g=\n_+\oplus\h\oplus\n_-$ and
set $\bb_\pm=\h\oplus\n_\pm$. Let $M(\lambda)$ be the Verma module
with the highest weight $\lambda\in\h^*$ and the highest weight
vector ${\mathbb I}_\lambda$. We call $\lambda$ {\it generic} if
$\langle\lambda, \alpha\rangle \notin {\mathbb Z}$ for any root
$\alpha$ of $\g$. It is well known that in this case $M(\lambda)$
is irreducible.

For any $\g$-module $V$ we will denote by $V[\mu]$ its subspace of
all weight vectors of weight $\mu\in\h^*$.

Now let us construct a map $\mbox{Hom}_\g(M(\lambda), M(\lambda)
\otimes F) \rightarrow F[0]$ for any $\lambda$. Choosing $\varphi
\in \mbox{Hom}_\g(M(\lambda), M(\lambda) \otimes F)$ we consider
$\varphi ({\mathbb I}_\lambda)\in M(\lambda) \otimes F$. Clearly,
$\varphi({\mathbb I}_\lambda) ={\mathbb I}_{\lambda} \otimes
f_\varphi + \sum_{\mu<\lambda} v_\mu \otimes f_\mu$, where
$v_\mu\in M(\lambda)[\mu]$. Obviously, $f_\varphi \in F[0]$. The
correspondence $\varphi\mapsto f_\varphi$ is the required map.

\begin{lemma}\label{l2}
If $M(\lambda)$ is irreducible, then this map is an isomorphism of
vector spaces.
\end{lemma}

\begin{proof}
Since $F$ is a direct sum of finite dimensional $\g$-modules, it
is enough to prove that $\mbox{Hom}_\g(M(\lambda),
M(\lambda)\otimes V) \cong V[0]$ if $\dim V<\infty$. This is well
known (see \cite{ES_DYBE}).
\end{proof}

\medskip

Lemma \ref{l2} provides a structure of associative algebra on
$F[0]$ since the space $\mbox{Hom}_\g(M(\lambda), M(\lambda)
\otimes F)$ has such a structure. We would like to describe this
structure in more details. Let us recall the definition of the
universal dynamical twist (see \cite{ES_DYBE}).

Let $\varphi \in \mbox{Hom}_\g(M(\lambda), M(\lambda) \otimes V)$
and $\psi \in \mbox{Hom}_\g(M(\lambda), M(\lambda) \otimes W)$,
where $V, W$ are finite dimensional $\g$-modules. We have the
following picture:
\begin{gather*}
M(\lambda) \stackrel{\psi}{\longrightarrow} M(\lambda) \otimes W
\stackrel{\varphi \otimes \id}{\longrightarrow} M(\lambda) \otimes
V\otimes W,\\
\psi({\mathbb I}_\lambda) = {\mathbb I}_\lambda \otimes u_\psi +
\sum_{\mu<\lambda} v_\mu\otimes u_{\mu,\psi},\\
\varphi({\mathbb I}_\lambda) = {\mathbb I}_\lambda \otimes
u_\varphi+\sum_{\mu<\lambda} v_\mu \otimes u_{\mu,\varphi},\\
(\varphi \otimes\id) \psi ({\mathbb I}_\lambda)={\mathbb
I}_\lambda \otimes u_{(\varphi \otimes
\id)\psi}+\sum_{\mu<\lambda} v_\mu \otimes u_{\mu,(\varphi\otimes
\id)\psi}.
\end{gather*}
It turns out that there exists a universal series $J(\lambda)\in
U\g\widehat{\otimes} U\g$ such that
\[
u_{(\varphi \otimes \id)\psi}=J(\lambda)_{V\otimes W}(u_\varphi
\otimes u_\psi).
\]
It is known that $J(\lambda) \in 1\otimes 1+(\n_-\cdot
U\n_-)\widehat{\otimes} (U\bb_+\cdot\n_+)$ and its coefficients
are rational functions on $\lambda$. Moreover, $J(\lambda)$ is a
dynamical twist.


\begin{proposition}
$f_{\varphi \ast
\psi}=(m\circ\overrightarrow{J(\lambda)})(f_\varphi \otimes
f_\psi)$.
\end{proposition}

\begin{proof} Taking into account that $F$ is a sum of finite
dimensional modules and the construction of $\varphi \ast \psi$ we
get the required result.
\end{proof}

\begin{corollary}\label{star_F_0}
The formula $f_1 \star_\lambda
f_2=(m\circ\overrightarrow{J(\lambda)})(f_1\otimes f_2)$ defines
an associative product on $F[0]$.\qed
\end{corollary}

In the sequel we denote by $F[0]_\lambda$ the obtained algebra.

\medskip

It is known that
\[
J({\lambda}/{\hbar})=1 \otimes 1+\hbar j(\lambda)+O(\hbar^2)
\]
and $r(\lambda)=j(\lambda)-j(\lambda)^{21}$ is the classical
triangular dynamical $r$-matrix (see \cite{ES_DYBE}). It was also
noticed in \cite{EKASrmatr} that $r(\lambda)$ defines a family of
$G$-invariant Poisson structures on $G/H$, which we denote by
$\{\cdot,\cdot\}_\lambda$. Any such a structure is in fact coming
from the Kirillov-Kostant-Souriau bracket on the coadjoint orbit
$O_\lambda \subset \g^\ast$.

\begin{corollary}
The multiplication $f_1 \star_{\lambda/\hbar}f_2$ is an
equivariant deformation quantization of the Poisson homogeneous
structure $\{\cdot,\cdot\}_\lambda$ on $G/H$ (and hence a
quantization of the Kirillov-Kostant-Souriau bracket on
$O_\lambda$).\qed
\end{corollary}


Now let us discuss the image of $U\g$ in $F[0]_\lambda$.

For this we have to compute $u_{x,\lambda}$ in
%
$\Ad_{g^{-1}} (x) {\mathbb I}_\lambda={\mathbb I}_\lambda \otimes
u_{x,\lambda}(g)+\ldots$.
%
We have the decomposition $U\g=U\h\oplus (\n_-\cdot
U\g+U\g\cdot\n_+)$, and for any $a \in U\g$ we define $(a)_0 \in
U\h$ as the corresponding projection. It is clear that
%
$u_{x,\lambda}(g)=(\Ad_{g^{-1}}(x))_0(\lambda)$
%
(we identify $U\h$ with polynomial functions on $\h^\ast$).

We get the following

\begin{proposition}
$u_{x,\lambda} \star_\lambda u_{y,\lambda}=u_{xy,\lambda}$.\qed
\end{proposition}

\begin{proposition}
The homomorphism $U\g \rightarrow F[0]_\lambda$ defined by $x
\mapsto u_{x,\lambda}$ is surjective for generic $\lambda$.
\end{proposition}

\begin{proof} We have the maps
\[
U\g\rightarrow U\g\otimes F \rightarrow U\h\otimes F[0]
\rightarrow F[0]
\]
defined by $x\mapsto \Ad_{g^{-1}}(x)\mapsto (\Ad_{g^{-1}}(x))_0
\mapsto (\Ad_{g^{-1}} (x))_0 (\lambda)$.

If we set $\deg F=0$, then the first two maps are filtered with
respect to the standard filtration of $U\g$.

If we go to the corresponding graded spaces we get the maps
\[
S(\g)\rightarrow S(\g) \otimes F\rightarrow S(\h)\otimes
F[0]\rightarrow F[0],
\]
$a \mapsto \Ad_{g^{-1}}(a) \mapsto (\Ad_{g^{-1}}(a))_0
\mapsto (\Ad_{g^{-1}}(a))_0 (\lambda)$. In this case the
surjectivity of the composition map $S(\g) \rightarrow F[0]$ for
generic $\lambda$ is well known.
It follows now that the map $U\g \rightarrow F[0]$ is also
surjective.
\end{proof}

\medskip

Let again $\lambda$ be generic. Since
$\Hom_\g(M(\lambda),M(\lambda)\otimes F)$ is isomorphic to
$F[0]_\lambda$ as an algebra, we get the following algebra
homomorphisms:
\[
U\g \rightarrow F[0]_\lambda \rightarrow \End M(\lambda),
\]
and the composition is the standard map $U\g \rightarrow \End
M(\lambda)$.

\begin{corollary}
For any generic $\lambda$ the images of $U\g$ and $F[0]_\lambda$
in $\End M(\lambda)$ coincide. \qed
\end{corollary}


\section{Verma modules and quantum homogeneous spaces}\label{q-Verma}

Now we are going to give an analogue of the results from Section
\ref{Verma} to the case of quantum universal enveloping algebras.
It will be convenient to develop first some formalism for a
general Hopf algebra (see Subsection \ref{Hopf}) and then proceed
to quantum universal enveloping algebras and corresponding Verma
modules (see Subsection \ref{QUE}).

\subsection{General construction}\label{Hopf}

Let $A$ be a Hopf algebra (over an arbitrary field $\Bbbk$). As
usual, we will denote by $\Delta$ (resp.\ $\varepsilon$, $S$) the
comultiplication (resp.\ counit, antipode) in $A$. We will
systematically use the Sweedler notation for comultiplication,
i.e., $\Delta(x)=\sum_{(x)}x_{(1)}\otimes x_{(2)}$,
$(\Delta\otimes\id)\Delta(x)=(\id\otimes\Delta)\Delta(x)=\sum_{(x)}x_{(1)}\otimes
x_{(2)}\otimes x_{(3)}$, etc.

Assume $M$ is a (left) $A$-module. We call an element $m\in M$
{\it locally finite} if $\dim Am<\infty$. Denote by $M_{\fin}$ the
subset of all locally finite elements in $M$. Clearly, $M_{\fin}$
is a submodule in $M$. Similarly, we can consider locally finite
elements in a right $A$-module $N$. For convenience, we will use
the notation $N^{\rrr}_{\fin}$ for the submodule of all locally
finite elements in this case.

Recall that the left (resp.\ right) adjoint action of $A$ on
itself is defined by the formula
$\ad_xa=\sum_{(x)}x_{(1)}aS(x_{(2)})$ (resp.\
$\ad^{\rrr}_xa=\sum_{(x)}S(x_{(1)})ax_{(2)}$). We denote by
$A_{\fin}$ (resp.\ $A^{\rrr}_{\fin}$) the corresponding submodules
of locally finite elements. Since
$\ad_x(ab)=\sum_{(x)}\ad_{x_{(1)}}(a)\ad_{x_{(2)}}(b)$, we see
that $A_{\fin}$ is a (unital) subalgebra in $A$; the same holds
for $A^{\rrr}_{\fin}$. If the antipode $S$ is invertible, then $S$
defines an isomorphism between $A_{\fin}$ and $A^{\rrr}_{\fin}$.
We will assume that $S$ is invertible.

Fix a Hopf subalgebra $F$ of the Hopf algebra $A^\star$ dual to
$A$.
In the sequel we
will use the left and right regular actions of $A$ on $F$ defined
respectively by the formulas $(\overrightarrow{a}f)(x)=f(xa)$ and
$(f\overleftarrow{a})(x)=f(ax)$.

Now let $M$ be a (left) $A$-module. Equip $F$ with the left
regular $A$-action and consider the space $\Hom_A(M,M\otimes F)$.
For any $\varphi,\psi\in\Hom_A(M,M\otimes F)$ define
$\varphi*\psi$ by formula \eqref{basic_star_prod}. It is
straightforward to verify that $\varphi*\psi\in\Hom_A(M,M\otimes
F)$, and this definition equips $\Hom_A(M,M\otimes F)$ with a
unital associative algebra structure.

Consider the map $\Phi:\Hom_A(M,M\otimes F)\to\End M$,
$\varphi\mapsto u_\varphi$, defined by
$u_\varphi(m)=(\id\otimes\varepsilon)(\varphi(m))$; here
$\varepsilon(f)=f(1)$ is the counit in $F$. In other words, if
$\varphi(m)=\sum_im_i\otimes f_i$, then
$u_\varphi(m)=\sum_if_i(1)m_i$. Using the fact that $\varepsilon$
is an algebra homomorphism it is easy to show that $\Phi$ is an
algebra homomorphism as well.

\begin{lemma}
The map $\Phi$ embeds $\Hom_A(M,M\otimes F)$ into $\End M$.
\end{lemma}

\begin{proof}
If $\varphi\in\Hom_A(M,M\otimes F)$, $\varphi(m)=\sum_im_i\otimes
f_i$, then
\[
\varphi(am)=a\varphi(m)=\sum_i\sum_{(a)}a_{(1)}m_i\otimes
\overrightarrow{a_{(2)}}f_i,
\]
and
\[
u_\varphi(am)=\sum_i\sum_{(a)}(\overrightarrow{a_{(2)}}f_i)(1)a_{(1)}m_i=
\sum_{(a)}a_{(1)}\left(\sum_if(a_{(2)})m_i\right).
\]
Assume now that $u_\varphi=0$, i.e.,
$\sum_{(a)}a_{(1)}\left(\sum_if(a_{(2)})m_i\right)=0$ for any
$a\in A$ and $m\in M$. Then, in particular,
\begin{gather*}
0=\sum_{(a)}S(a_{(1)})a_{(2)}\left(\sum_if(a_{(3)})m_i\right)=\\
\sum_{(a)}\varepsilon(a_{(1)})\left(\sum_if(a_{(2)})m_i\right)=\sum_if_i(a)m_i
\end{gather*}
for any $a\in A$ and $m\in M$. Obviously, this means that
$\varphi=0$.
\end{proof}

From now on we assume that $F$ contains all matrix elements of the
(left) adjoint action of $A$ on $A_{\fin}$. Since $F$ is closed
under the antipode $(Sf)(x)=f(S(x))$, we see that this assumption
is equivalent to the fact that $F$ contains all matrix elements of
the right adjoint action of $A$ on $A^{\rrr}_{\fin}$.

Let $a\in A^{\rrr}_{\fin}$, i.e., for any $x\in A$ we have
$\ad^{\rrr}_xa=\sum_if_i(x)a_i$, where $f_i\in A^\star$, $a_i\in
A$. In fact, we see that $f_i\in F$ by the assumption above.
Define a linear map $\varphi_a:M\to M\otimes F$ by the formula
$\varphi_a(m)=\sum_ia_im\otimes f_i$. Clearly, $\varphi_a$ is well
defined.

\begin{lemma}
For any $a\in A^{\rrr}_{\fin}$ we have
$\varphi_a\in\Hom_A(M,M\otimes F)$.
\end{lemma}

\begin{proof}
Let $b\in A$. Notice that
\[
\sum_{(b)}b_{(1)}\ad^{\rrr}_{b_{(2)}}y=\sum_{(b)}b_{(1)}S(b_{(2)})yb_{(3)}=y\sum_{(b)}\varepsilon(b_{(1)})b_{(2)}=yb
\]
for any $y\in A$. Therefore for any $x\in A$ we have
\begin{gather*}
\sum_if_i(x)a_ib=(\ad^{\rrr}_xa)b=\sum_{(b)}b_{(1)}\ad^{\rrr}_{b_{(2)}}\ad^{\rrr}_xa=
\sum_{(b)}b_{(1)}\ad^{\rrr}_{xb_{(2)}}a=\\
\sum_{(b)}b_{(1)}\left(\sum_if_i(xb_{(2)})a_i\right)=
\sum_{(b)}\sum_i(\overrightarrow{b_{(2)}}f_i)(x)b_{(1)}a_i,
\end{gather*}
and
\[
\varphi_a(bm)=\sum_ia_ibm\otimes f_i=
\sum_{(b)}\sum_ib_{(1)}a_im\otimes\overrightarrow{b_{(2)}}f_i=b\varphi_a(m).
\]
\end{proof}

Denote by $\Psi:A^{\rrr}_{\fin}\to\Hom_A(M,M\otimes F)$ the linear
map constructed above (i.e., $\Psi:a\mapsto\varphi_a$).

\begin{lemma}
The map $\Psi$ is an algebra homomorphism.
\end{lemma}

\begin{proof}
Let $a,b\in A^{\rrr}_{\fin}$, $x\in A$,
$\ad^{\rrr}_xa=\sum_if_i(x)a_i$, $\ad^{\rrr}_xb=\sum_jg_j(x)b_j$.
Then
\begin{gather*}
\ad^{\rrr}_x(ab)=\sum_{(x)}\ad_{x_{(1)}}(a)\ad_{x_{(2)}}(b)=\\
\sum_{i,j}\sum_{(x)}f_i(x_{(1)})g_j(x_{(2)})a_ib_j=\sum_{i,j}(f_ig_j)(x)a_ib_j.
\end{gather*}
Thus
\[
\varphi_{ab}(m)=\sum_{i,j}a_ib_jm\otimes f_ig_j=
(\varphi_a*\varphi_b)(m)
\]
for any $m\in M$.
\end{proof}

\begin{remark}
It follows directly from the definitions that the composition
$\Phi\Psi$ equals the restriction to $A^{\rrr}_{\fin}$ of the
canonical homomorphism $A\to\End M$, $a\mapsto a_M$.
\end{remark}

Now consider $A^{\rrr}_{\fin}$, $\Hom_A(M,M\otimes F)$ and $\End
M$ as right $A$-modules: $A^{\rrr}_{\fin}$ via right adjoint
action, $\Hom_A(M,M\otimes F)$ via right regular action on $F$
(i.e., $(\varphi\cdot a)(m)=(\id\otimes
\overleftarrow{a})(\varphi(m))$), and $\End M$ in a standard way
(i.e., $u\cdot a=\sum_{(a)}S(a_{(1)})_Mu{a_{(2)}}_M$). Note that
$A^{\rrr}_{\fin}$, $\Hom_A(M,M\otimes F)$ and $\End M$ equipped
with these structures are indeed right $A$-module algebras, i.e.,
the multiplication map is a module morphism, and the unit is
invariant.

\begin{lemma}
The maps $\Phi$ and $\Psi$ are morphisms of right $A$-modules.
\end{lemma}

\begin{proof}
Straightforward.
\end{proof}

\begin{corollary}\label{Phi_Psi_q}
We have the following morphisms of right $A$-module algebras:
\[
A^{\rrr}_{\fin}\stackrel{\Psi}{\longrightarrow}\Hom_A(M,M\otimes
F)^{\rrr}_{\fin}\stackrel{\Phi}{\longrightarrow}(\End
M)^{\rrr}_{\fin},
\]
and $\Phi\Psi$ is the restriction of the canonical morphism $A\to
\End M$. \qed
\end{corollary}

\subsection{QUE algebra case}\label{QUE}

Now suppose $A=\check{U}_q\g$,
where $\g$ is a complex simple Lie algebra (see \cite[\S
3.2.10]{Joseph_book} or \cite{Joseph_Letzter_separation}). We
consider $A$ as an algebra over $\overline{\mathbb C(q)}$, the
algebraic closure of the field $\mathbb C(q)$ of rational
functions on the indeterminate $q$. It is known that the adjoint
action of $A$ is not locally finite. The subalgebra
$A_{\fin}\subset A$ was studied in \cite{Joseph_Letzter_adjoint,
Joseph_Letzter_separation} (see also \cite{Joseph_book}).

Let $F=\mathbb C[G]_q$ be the quantized algebra of regular
functions on an algebraic group $G$ corresponding to $\g$ (see
\cite{Joseph_book, KS_alg_func}). We can consider $F$ as a Hopf
subalgebra in $A^\star$. Clearly, $F$ satisfies the requirements
of the previous subsection. Notice that $F$ is a sum of finite
dimensional admissible $A$-modules with respect to both left and
right regular actions of $A$ (see \cite{KS_alg_func}).

Denote by $\h^*_{\mathbb Q}$ the $\mathbb Q$-span of the weight
lattice of $\g$. Consider the Verma module $M(\lambda)$ for $A$
with the highest weight $\lambda\in\h^*_{\mathbb Q}$ and the
highest weight vector ${\mathbb I}_\lambda$. As in the classical
case, if $\lambda$ is generic (i.e., $\langle\lambda,
\alpha\rangle \notin {\mathbb Z}$ for any root $\alpha$ of $\g$),
then $M(\lambda)$ is irreducible.

For any left $A$-module $V$ we will denote by $V[\mu]$ its
subspace of all weight vectors of weight $\mu\in\h^*_{\mathbb Q}$.

Now for any $\lambda\in\h^*_{\mathbb Q}$ we construct a linear map
\[
\Theta_\lambda:\Hom_A(M(\lambda), M(\lambda)\otimes F)\to F[0]
\]
(here $F$ is considered as an $A$-module via left regular action).
Take $\varphi \in \Hom_A(M(\lambda), M(\lambda) \otimes F)$ and
consider $\varphi ({\mathbb I}_\lambda)\in M(\lambda) \otimes F$.
Clearly, $\varphi({\mathbb I}_\lambda) ={\mathbb I}_{\lambda}
\otimes f_\varphi + \sum_{\mu<\lambda} v_\mu \otimes f_\mu$, where
$v_\mu\in M(\lambda)[\mu]$. We see that $f_\varphi \in F[0]$. The
correspondence $\varphi\mapsto f_\varphi$ is the map of concern.
Notice that both $\Hom_A(M(\lambda), M(\lambda)\otimes F)$ and
$F[0]$ are right $A$-modules (via right regular action of $A$ on
$F$), and $\Theta_\lambda$ is compatible with these structures.

By the same argument as in the classical case (see Lemma
\ref{l2}), we have the following

\begin{lemma}\label{l2-q}
If $M(\lambda)$ is irreducible, then $\Theta_\lambda$ is an
isomorphism (of right $A$-modules).\qed
\end{lemma}

Now for any generic $\lambda$ we can use $\Theta_\lambda$ to
transfer to $F[0]$ the product $*$ on $\Hom_A(M(\lambda),
M(\lambda)\otimes F)$, which was constructed in the previous
subsection.

Arguing like in the classical case, we see that $f_{\varphi \ast
\psi}=(m\circ\overrightarrow{J_q(\lambda)})(f_\varphi \otimes
f_\psi)$ for any $\varphi,\psi\in\Hom_A(M(\lambda),
M(\lambda)\otimes F)$. Here $J_q(\lambda)$ is a universal quantum
dynamical twist for $A$ (see \cite{ES_DYBE}). Therefore we get

\begin{corollary}\label{star_F_0_q}
The formula $f_1 \star_\lambda
f_2=(m\circ\overrightarrow{J_q(\lambda)})(f_1\otimes f_2)$ defines
an associative product on $F[0]$, and
$\Theta_\lambda:(\Hom_A(M(\lambda), M(\lambda)\otimes F),*)\to
(F[0],\star_\lambda)$ is an algebra isomorphism.\qed
\end{corollary}

Let us denote by $F[0]_\lambda$ the algebra
$(F[0],\star_\lambda)$.

\begin{remark}
Obviously, $F[0]_\lambda$ is a right $A$-module algebra (i.e.,
$(f_1 \star_\lambda
f_2)\overleftarrow{a}=\sum_{(a)}f_1\overleftarrow{a_{(1)}}\star_\lambda
f_2\overleftarrow{a_{(2)}}$ and
$1\overleftarrow{a}=\varepsilon(a)1$ for any $a\in A$).
\end{remark}

Let us identify $\Hom_A(M(\lambda), M(\lambda)\otimes F)$ and
$F[0]_\lambda$ via $\Theta_\lambda$. Note that the right regular
action on $F[0]\subset F$ is locally finite, i.e.,
$\left(F[0]\right)^{\rrr}_{\fin}=F[0]$. By Corollary
\ref{Phi_Psi_q} we have
\[
A^{\rrr}_{\fin}\stackrel{\Psi}{\longrightarrow}F[0]_\lambda\stackrel{\Phi}{\longrightarrow}(\End
M(\lambda))^{\rrr}_{\fin},
\]
and $\Phi\Psi$ is the restriction of the canonical map $A\to\End
M(\lambda)$. It is known that this restriction is surjective (cf.\
\cite{Joseph_book, Joseph_Letzter_Verma}). Since $\Phi$ is an
embedding, we see that the following holds:

\begin{theorem}
The map $\Phi$ defines an isomorphism between right $A$-module
algebras $F[0]_\lambda$ and $(\End M(\lambda))^{\rrr}_{\fin}$.\qed
\end{theorem}

\begin{remark}
Let $G$ be a Lie group corresponding to $\g$, and $H\subset G$ a
Cartan subgroup. For each $\lambda$ the right $A$-module algebra
$F[0]_\lambda$ is a ``quantization'' of a $(G,\pi_0)$-homogeneous
Poisson structure on $G/H$, where $\pi_0$ is the Poisson Lie group
structure on $G$ defined by the standard quasitriangular solution
of the classical Yang-Baxter equation for $\g$.
\end{remark}

\section{From dynamical to non-dynamical twists}\label{ndtwists}

Let $\g$ be a Lie algebra, $\h$ its abelian subalgebra. Suppose
\begin{equation*}
J:\h^*\to (U\g\otimes U\g)^\h[[\hbar]]
\end{equation*}
is a quantum dynamical twist.

Assume that there exists a subalgebra $\v\subset\g$ such that
$\g=\h\oplus\v$ as a vector space. Notice that $U\g=U\v\oplus
U\g\cdot\h$.
Denote by $J_\v(\lambda)$
the image of $J(\lambda)$ under the projection onto $(U\v\otimes
U\v)[[\hbar]]$ along $\left((U\g\cdot\h\otimes U\g)\oplus
(U\g\otimes U\g\cdot\h)\right)[[\hbar]]$.

\begin{theorem}
For any $\lambda\in\Dom J$ the element $J_\v(\lambda)$ is a
quantum twist for $U\v$, i.e.,
$J_\v(\lambda)^{12,3}J_\v(\lambda)^{12}=J_\v(\lambda)^{1,23}J_\v(\lambda)^{23}$,
and
$(\varepsilon\otimes\id)(J_\v(\lambda))=(\id\otimes\varepsilon)(J_\v(\lambda))=1$.
\end{theorem}

\begin{proof}
We have $J(\lambda)=J_\v(\lambda)+J_\h(\lambda)$, where
\begin{equation*}
J_\h(\lambda)\in\left((U\g\cdot\h\otimes U\g)\oplus (U\g\otimes
U\g\cdot\h)\right)[[\hbar]].
\end{equation*}
Denote by $A_\v$ the projection of $A\in(U\g\otimes U\g\otimes
U\g)[[\hbar]]$ onto $(U\v\otimes U\v\otimes U\v)[[\hbar]]$ along
\begin{equation*}
\left((U\g\cdot\h\otimes U\g\otimes U\g)\oplus(U\g\otimes
U\g\cdot\h\otimes U\g)\oplus(U\g\otimes U\g\otimes
U\g\cdot\h)\right)[[\hbar]].
\end{equation*}
One can calculate
directly, using the fact that $U\v$ (resp.\ $U\g\cdot\h$) is a
subalgebra and a coideal (resp.\ a left ideal and a coideal) in
$U\g$, that
\begin{equation*}
(J(\lambda)^{12,3}J(\lambda-\hbar
h^{(3)})^{12})_\v=J_\v(\lambda)^{12,3}J_\v(\lambda)^{12}+(J_\h(\lambda)^{12,3}J_\v(\lambda)^{12})_\v
\end{equation*}
and
\begin{equation*}
(J(\lambda)^{1,23}J(\lambda)^{23})_\v=
J_\v(\lambda)^{1,23}J_\v(\lambda)^{23}+(J_\h(\lambda)^{1,23}J_\v(\lambda)^{23})_\v.
\end{equation*}
Therefore from (\ref{qdtwist}) it follows that
\begin{gather*}
J_\v(\lambda)^{12,3}J_\v(\lambda)^{12}+(J_\h(\lambda)^{12,3}J_\v(\lambda)^{12})_\v=\\
J_\v(\lambda)^{1,23}J_\v(\lambda)^{23}+(J_\h(\lambda)^{1,23}J_\v(\lambda)^{23})_\v.
\end{gather*}

Let us prove that in fact
$(J_\h(\lambda)^{12,3}J_\v(\lambda)^{12})_\v=
(J_\h(\lambda)^{1,23}J_\v(\lambda)^{23})_\v=0$.

Let
\[
J_\h(\lambda)=\sum_{m\geq0}J_\h^{(m)}(\lambda)\hbar^m,\ \
J_\v(\lambda)=\sum_{n\geq0}J_\v^{(n)}(\lambda)\hbar^n,
\]
where
$J_\h^{(m)}(\lambda)\in(U\g\cdot\h\otimes U\g)\oplus (U\g\otimes
U\g\cdot\h)$, $J_\v^{(n)}(\lambda)\in U\v\otimes U\v$.

Clearly, it is enough to show that for any $m,n$ we have
\[
\left(J_\h^{(m)}(\lambda)^{12,3}J_\v^{(n)}(\lambda)^{12}\right)_\v=
\left(J_\h^{(m)}(\lambda)^{1,23}J_\v^{(n)}(\lambda)^{23}\right)_\v=0.
\]

Indeed, write
\begin{equation*}
J_\v^{(n)}(\lambda)=\sum_ix_i\otimes y_i,
\end{equation*}
where
$x_i,y_i\in U\v$, and
\begin{equation*}
J_\h^{(m)}(\lambda)=\sum_j(a_jh_j'\otimes b_j+c_j\otimes
d_jh_j''),
\end{equation*}
where $a_j,b_j,c_j,d_j\in U\g$, $h_j',h_j''\in\h$. We have
\begin{gather*}
J_\h^{(m)}(\lambda)^{12,3}=(\Delta\otimes\id)(J_\h^{(m)}(\lambda))=\\
\sum_j\left((\Delta(a_j)(h_j'\otimes1+1\otimes
h_j'))\otimes b_j+\Delta(c_j)\otimes d_jh_j''\right)
\end{gather*}
and
\begin{equation*}
J_\v^{(n)}(\lambda)^{12}=J_\v^{(n)}(\lambda)\otimes1=\sum_ix_i\otimes
y_i\otimes1.
\end{equation*}
Therefore
\begin{gather*}
(J_\h^{(m)}(\lambda)^{12,3}J_\v^{(n)}(\lambda)^{12})_\v=\\
\left(\sum_{i,j}\left((\Delta(a_j)(h_j'\otimes1+1\otimes
h_j'))\otimes b_j\right)\left(x_i\otimes
y_i\otimes1\right)\right)_\v=\\
\left(\sum_{i,j}\left(\Delta(a_j)(h_j'x_i\otimes y_i+x_i\otimes
h_j'y_i)\right)\otimes b_j\right)_\v=\\
\left(\sum_{i,j}\left(\Delta(a_j)([h_j',x_i]\otimes y_i+x_i\otimes
[h_j',y_i])\right)\otimes b_j\right)_\v=\\
\left(\sum_{i,j}\left(\Delta(a_j)([h_j',x_i]_\v\otimes
y_i+x_i\otimes
[h_j',y_i]_\v)\right)\otimes b_j\right)_\v=\\
\left(\sum_j\Delta(a_j)\otimes
b_j\cdot\left(\left(\sum_i\left([h_j',x_i]_\v\otimes
y_i+x_i\otimes
[h_j',y_i]_\v\right)\right)\otimes1\right)\right)_\v;
\end{gather*}
here $[h_j',x_i]_\v$ means the projection of $[h_j',x_i]$ onto
$U\v$ along $U\g\cdot\h$, etc.

Now recall that $\ad_h(J(\lambda))=0$ for all $h\in\h$. Projecting
this equation onto $(U\v\otimes U\v)[[\hbar]]$ along
$\left((U\g\cdot\h\otimes U\g)\oplus (U\g\otimes
U\g\cdot\h)\right)[[\hbar]]$, we get
\begin{equation*}
\sum_i\left([h,x_i]_\v\otimes y_i+x_i\otimes [h,y_i]_\v\right)=0
\end{equation*}
for all $h\in\h$. Combining this with the previous computation, we
see that
$(J_\h^{(m)}(\lambda)^{12,3}J_\v^{(n)}(\lambda)^{12})_\v=0$.

Similarly,
$(J_\h^{(m)}(\lambda)^{1,23}J_\v^{(n)}(\lambda)^{23})_\v=0$.

Finally, the counit condition on $J_\v(\lambda)$ follows easily
from (\ref{qdtwist2}).
\end{proof}

\begin{example}
Suppose $\g=\sl(2)$. Let $x,y,h$ be the standard basis in $\g$,
and $\h=\mathbb Ch$. Notice that $\g=\h\oplus\v$, where
\begin{equation*}
\v=
g\left\{\left(\begin{array}{cc}*&*\\0&*\end{array}\right)\right\}g^{-1},\
\ g=\left(\begin{array}{cc}1&0\\1&1\end{array}\right).
\end{equation*}

Consider the ABRR quantum dynamical twist $J$ for $(\g,\h)$ (see
Example \ref{sl2ex}). We have $y=b+c$, $x=a-c$, where
$b=y+\frac{1}{2}h=g\left(\frac{1}{2}h\right)g^{-1}\in\v$,
$a=x-\frac{1}{2}h=g\left(x+\frac{1}{2}h\right)g^{-1}\in\v$,
$c=-\frac{1}{2}h\in\h$. Obviously, $[c,b]=b+c$ and $[-c,a]=a-c$.

\begin{lemma}
The projection of $y^n=(b+c)^n$ onto $U\v$ along $U\g\cdot\h$
equals $b(b+1)\ldots(b+n-1)$.
\end{lemma}

\begin{proof}
One can easily verify by induction that
\begin{equation*}
cb^n=b\left((b+1)^n-b^n\right)+(b+1)^nc.
\end{equation*}
Therefore
\begin{equation*}
cf(b)=b\left(f(b+1)-f(b)\right)+f(b+1)c
\end{equation*}
for any polynomial $f$. Finally,
\begin{gather*}
(b+c)\cdot b(b+1)\ldots(b+n-1)=\\
b^2(b+1)\ldots(b+n-1)+\\b\left((b+1)(b+2)\ldots(b+n)-b(b+1)\ldots(b+n-1)\right)+\\(b+1)(b+2)\ldots(b+n)c=\\
b(b+1)(b+2)\ldots(b+n)+(b+1)(b+2)\ldots(b+n)c.
\end{gather*}
\end{proof}

Applying the lemma, we see that
\begin{equation*}
J_\v(\lambda)=1+\sum_{n\geq1}
\frac{(-1)^n\hbar^nv_n}{n!\lambda(\lambda-\hbar)\ldots(\lambda-(n-1)\hbar)},
\end{equation*}
where
\begin{gather*}
v_n=b(b+1)\ldots(b+n-1) \otimes a(a+1)\ldots(a+n-1).
\end{gather*}
\end{example}

\end{sloppy}

\small

\noindent E.K.: Department of Mathematics, Kharkov National University,\\
4 Svobody Sq., Kharkov \,61077, Ukraine;\\
Institute for Low Temperature Physics \& Engineering,\\
47 Lenin Avenue, Kharkov \,61103, Ukraine\\
e-mail: {\small \tt eugene.a.karolinsky@univer.kharkov.ua}

\medskip

\noindent K.M.: Department of Mathematics, Kharkov National University,\\
4 Svobody Sq., Kharkov \,61077, Ukraine

\medskip

\noindent A.S.: Department of Mathematics, University of G\"oteborg,\\
SE-412 96 G\"oteborg, Sweden\\
e-mail: {\small \tt astolin@math.chalmers.se}

\medskip

\noindent V.T.: St.\,Petersburg Branch of Steklov Mathematical Institute,\\
Fontanka 27, St.\,Petersburg \,191023, Russia;\\
Department of Mathematical Sciences, IUPUI,\\
Indianapolis, IN 46202, USA \\
e-mail: {\small \tt vt@pdmi.ras.ru; vt@math.iupui.edu}

\end{document}